\documentclass[12pt]{article}
\title{Limits of Thompson's group F}
\author{Roland Zarzycki}
\usepackage[latin2]{inputenc}
\usepackage{amssymb}
\usepackage{indentfirst}
\usepackage{latexsym}
\begin{document}

\maketitle

\newtheorem{theo}{Theorem}[section]
\newtheorem{lemm}[theo]{Lemma}
\newtheorem{coro}[theo]{Corollary}
\newtheorem{defi}[theo]{Definition}
\newtheorem{prop}[theo]{Proposition}
\newtheorem{fact}[theo]{Fact}
\newtheorem{rema}[theo]{Remark}
\newtheorem{exam}[theo]{Example}
\newtheorem{ques}[theo]{Question}
\newtheorem{clai}[theo]{Claim}

\begin{abstract}

 Let $F$ be the Thompson's group $\langle x_{0}, x_{1}| [x_{0}x_{1}^{-1}, x_{0}^{-i}x_{1}
x_{0}^{i}], i=1,2 \rangle$. Let $G_{n}= \langle y_{1},\ldots , y_{m}, x_{0},x_{1} \vert [x_{0}
x_{1}^{-1},x_{0}^{-i}x_{1}x_{0}^{i}], y_{j}^{-1}$ $g_{j,n}(x_{0},x_{1}), i=1, 2,$ $j\leq m \rangle$,
where $g_{j,n}(x_{0},x_{1})\in F$, $n\in\mathbb{N}$, be a family of groups isomorphic to
$F$ and marked by $m+2$ elements. If the sequence $(G_{n})_{n<\omega}$ is
convergent in the space of marked groups and $G$ is the corresponding limit we say
that $G$ is an $F$-limit group. The paper is devoted to a description of $F$-limit groups.

\end{abstract}

\section{Preliminaries}

 The notion of \emph{limit group} was introduced by Z. Sela in his work on characterization
of elementary equivalence of free groups \cite{Sel}. This approach has been extended in the
paper of C. Champetier and V. Guirardel \cite{GC}, where the authors look at \emph{limit
groups} as limits of convergent sequences in a space of \emph{marked groups}. They have
given a description of Sela's limit groups in these terms (with respect to the class of free
groups). This approach has been also aplied by L. Guyot and Y. Stalder \cite{SG} to the
class of Baumslag-Solitar groups. \parskip0pt

 Thompson's group $F$ has remained one of the most interesting objects in geometric group
theory. We study $F$-limit groups. We show in this paper, that among $F$-limit groups there
are no free products of $F$ with any non-trivial group. Moreover, we prove that among
$F$-limit groups there are no HNN-extensions over cyclic subgroups. \parskip0pt

 In the remaining part of the section we recollect some useful definitions and facts
concerning limit groups and Thompson's group $F$. In Section 2 we present results
concerning free products and in Section 3 results concerning HNN-extensions. \parskip0pt

 A \emph{marked group} $(G,S)$ is a group G with a distinguished set of generators $S =
(s_{1}, s_{2}, \ldots , s_{n})$. For fixed $n$, let $\mathcal{G} _{n}$ be the set of all
$n$-generated groups marked by $n$ generators (up to isomorphism of marked groups).
Following \cite{GC} we put certain metric on $\mathcal{G} _{n}$. We will say that two marked
groups $(G,S), (G',S')\in\mathcal{G} _{n}$ are at distance less or equal to $e^{-R}$ if they
have exactly the same relations of length at most $R$. The set $\mathcal{G} _{n}$ equiped
with this metric is a compact space \cite{GC}. \emph{Limit groups} are simply limits of
convergent sequences in this metric space.

\begin{defi} Let $G$ be an $n$-generated group. A marked group in $\mathcal{G} _{n}$ is
	a $G$-limit group if it is a limit of marked groups each isomorphic to $G$.
\end{defi}

 To introduce the Thompson's group $F$ we will follow \cite{CFP}.

\begin{defi} Thompson's group $F$ is the group given by the following infinite group
presentation:
	$$ \langle x_{0}, x_{1}, x_{2}, \ldots | x_{j}x_{i}=x_{i}x_{j+1} (i<j) \rangle $$
\end{defi}

 In fact $F$ is finitely presented:
	$$ F = \langle x_{0}, x_{1}| [x_{0}x_{1}^{-1}, x_{0}^{-i}x_{1}x_{0}^{i}], i=1,2
\rangle . $$
 Every non-trivial element of $F$ can be uniquely expressed in the normal form:
	$$ x_{0}^{b_{0}}x_{1}^{b_{1}}x_{2}^{b_{2}}\ldots x_{n}^{b_{n}}x_{n}^{-a_{n}}\ldots
	x_{2}^{-a_{2}}x_{1}^{-a_{1}}x_{0}^{-a_{0}},$$
where $n$, $a_{0}, \ldots , a_{n}$, $b_{0}, \ldots , b_{n}$ are non-negative integers such
that: \\
	i) exactly one of $a_{n}$ and $b_{n}$ is nonzero; \\
	ii) if $a_{k}>0$ and $b_{k}>0$ for some integer $k$ with $0\leq k < n$, then
	$a_{k+1}>0$ or $b_{k+1}>0$. \parskip0pt

 We study properties of $F$-limit groups. For this purpose let us consider a sequence,
$(g_{i,n})_{n<\omega}$, $1\leq i\leq t$, of elements taken from the group $F$ and the
corresponding sequence of limit groups marked by $t+2$ elements, $G_{n} = (F,(x_{0},x
_{1},g_{1,n},\ldots ,g_{t,n}))$, $n\in\mathbb{N}$, where $x_{0}$ and $x_{1}$ are the
standard generators of $F$. Assuming that such a sequence is convergent in the space
of groups marked by $t+2$ elements, denote by $G = (\langle x_{0}, x_{1}, g_{1}\ldots ,
g_{t}|R_{F} \cup R_{G} \rangle, (x_{0}, x_{1}, g_{1},\ldots , g_{t}))$ the limit group formed
in that manner; here $x_{0}$, $x_{1}$ are "limits" of constant sequences $(x_{0})_{n<
\omega}$ and $(x_{0})_{n<\omega}$, $g_{i}$ is the "limit" of $(g_{i,n})_{n<\omega}$ for
$1\leq i\leq t$, $R_{F}$ and $R_{G}$ refer respectively to the set of standard relations
taken from $F$ and the set (possibly infinite) of new relations. \parskip0pt

 It has been shown in \cite{GC} that in the case of free groups some standard constructions
can be obtained as limits of free groups. For example, it is possible to get $\mathbb{Z}
^{k}$ as a limit of $\mathbb{Z}$ and $\mathbb{F}_{k}$ as a limit of $\mathbb{F} _{2}$.  On
the other hand, the direct product of $\mathbb{F} _{2}$ and $\mathbb{Z}$ can not be
obtained as a limit group. HNN-extensions often occure in the class of limit groups (with
respect to free groups). For example, the following groups are the limits of convergent
sequences in the space of free groups marked by three elements: the free group of rank
$3$, the free abelian group of rank $3$ or a HNN-extension over a cyclic subgroup of the
free group of rank $2$ (\cite{FGM}). All non-exceptional surface groups form another broad
class of interesting examples (\cite{BaB}, \cite{BaG}). \parskip0pt

 In the case of Thompson's group the situation is not so clear. Since the centrum of $F$ is
trivial it is surely not possible to obtain any direct product with the whole group as an $F$-limit
group. In 1985 Brin and Squier \cite{BS} showed that Thompson's group $F$ does not satisfy
any law (also see Abert's paper \cite{A} for a shorter proof). However, in this paper we show
that there are certain non-trivial words with constants over $F$ (which will be called later laws
with constants), which are equal to the identity for each evaluation in $F$. This implies that no
free product of $F$ with any non-trivial group is admissible as limit group with respect to $F$
(see Section 2). Moreover, we prove that HNN-extensions over a cyclic subgroup are not
admissible as limit groups with respect to $F$ (see Section 3). \parskip0pt

 There are many geometric interpretations of $F$, but here we will use the following one.
Consider the set of all strictly increasing continuous piecewise-linear functions from the
closed unit interval onto itself. Then the group $F$ is realized by the set of all such
functions, which are differentiable except at finitely many dyadic rational numbers and such
that all slopes (deriviatives) are integer powers of 2. The corresponding group operation is
just the composition. For the further reference it will be usefull to give an explicit form of
the generators $x_{0}, x_{1}, \ldots$ in terms of piecewise-linear functions:

$$ x_{n}(t) = \left\{ \begin{array}{ll} t & \textrm{, $t\in [0,\frac{2^{n}-1}{2^{n}} ]$} \\
	\frac{1}{2}t + \frac{2^{n}-1}{2^{n+1}} & \textrm{, $t\in [\frac{2^{n}-1}{2^{n}}, \frac{2^{n+1}-1}
	{2^{n+1}} ]$} \\ t - \frac{1}{2^{n+2}} & \textrm{, $t\in [\frac{2^{n+1} -1}{2^{n+1}},
	\frac{2^{n+2}-1}{2^{n+2}}]$} \\	2t-1 & \textrm{, $t\in [\frac{2^
	{n+2}-1}{2^{n+2}},1]$} \end{array}\right. $$ for $n = 0, 1,\ldots$. \parskip0pt

 For any diadic subinterval $[a,b]\subset [0,1]$, let us consider the set, of elements in $F$,
which are trivial on its complement, and denote it by $F _{[a,b]}$. We know that it forms a
subgroup of $F$, which is isomorphic to the whole group. Let us denote its standard infinite
set of generators by $x_{[a,b],0}, x_{[a,b],1}, x_{[a,b],2}, \ldots$.

 Let us consider an arbitrary element $g$ in $F$ and treat it as a piecewise-linear
homeomorphism of the interval $[0,1]$. Let $supp(g)$ be the set $\{ x\in [0,1] : g(x) \neq x \}$
and $\overline{supp}(g)$ the topological closure of $supp(g)$. We will call each point from the
set $P_{g} = (\overline{supp}(g)\setminus supp(g))\cap\mathbb{Z} [\frac{1}{2}]$ a \emph{dividing
point} of $g$. This set is obviously finite and thus we get a finite subdivision of $[0,1]$ of the
form $[0=p_{0}, p_{1}], [p_{1}, p_{2}], \ldots , [p_{n-1}, p_{n}=1]$ for some natural $n$. It is
easy to see that $g$ can be presented as $g = g_{1}g_{2}\ldots g_{n}$, where $g_{i}\in F
_{[p_{i-1}, p_{i}]}$ for each $i$. Since $g$ can act trivially on some of these subintervals, some
of the elements $g_{1},\ldots ,g_{n}$ may be trivial. We call the set of all non-trivial elements
from $\{ g_{1},\ldots ,g_{n}\}$ the \emph{defragmentation} of $g$.

\begin{fact}[Corollary 15.36 in \cite{GS}, Proposition 3.2 in \cite{KM}]\label{GS}
	The centralizer of any element $g\in F$ is the direct product of finitely many cyclic
	groups and finitely many groups isomorphic to $F$.
\end{fact}

 Moreover if the element $g\in F$ has the defragmentation $g=g_1 \ldots g_n$, then some roots
of the elements $g_1, \ldots , g_n$ are the generators of cyclic components of the decomposition
of the centralizer above. The components of this decomposition which are isomorphic to $F$ are
just the groups of the form $F_{[a,b]}$, where $[a,b]$ is one of the subintervals $[p_{i-1},p_{i}]
\subset [0,1]$, which are stabilized pointwise by $g$. Generally, if we interpret the elements of
$F$ as functions, the relations occuring in the presentation of $F$, $[x_{0}x_{1}^{-1}, x_{0}^{-i}
x_{1}x_{0}^{i}]$ for $i = 1,2$, have to assure, that two functions, which have mutually disjoint
supports except of finitely many points, commute. In particular, these relations imply analogous
relations for different $i>2$. According to the fact that $x_{0}^{-i}x_{1}x_{0}^{i} = x_{i+1}$,
we conclude that all the relations of the form $[x_{0}x_{1}^{-1}, x_{M}]$, $M>1$, hold in
Thompson's group $F$. We often refer to these geometrical observations.

 I am grateful to the referee for his helpful remarks.

\section{Free products}

 Brin and Squier have shown in $\cite{BS}$ that the Thompson's group $F$ does not satisfy any
group law. In this section we show how to construct words with constants from $F$, which are equal
to the identity for any substitution in $F$.

\begin{defi}
	 Let $w(y_{1},\ldots , y_{t})$ be a non-trivial word over $F$, reduced in the group
	$\mathbb{F} _{t}\ast F$ and containing at least one variable. We will call $w$ a
	\emph{law with constants} in $F$ if for any $\bar{g}=(g_{1},\ldots , g_{t})\in F^{t}$,
	the value $w(\bar{g})$ is equal to $1_{F}$.
\end{defi}

 The following proposition gives a construction of certain laws with constants in $F$.

\begin{prop} \label{lwc}
	 Consider the standard action of Thompson's group $F$ on $[0,1]$. Suppose we are given
	four pairwise disjoint closed diadic subintervals $I_{i}=[p_{i}, q_{i}]\subset [0,1]$, $1\leq i
	\leq 4$, and assume that $p_{1}<p_{2}<p_{3}<p_{4}$. Then for any non-trivial $h_{1}\in F
	_{I_{1}}$, $h_{2}\in F_{I_{2}}$, $h_{3}\in F_{I_{3}}$ and $h_{4}\in F_{I_{4}}$, the word $w$
	obtained from
	$$[y^{-1}h_{1}^{-1}yh_{4}^{-1}y^{-1}h_{1}yh_{4}, y^{-1}h_{2}^{-1}yh_{3}^{-1}y^{-1}h_{2}yh_{3}]$$
	by reduction in $\mathbb{Z}\ast F$ (we treat the variable $y$ as a generator of
	$\mathbb{Z}$) is a law with constants in $F$.
\end{prop}

\emph{Proof.} \ We will use the following notation: $w_{14}=y^{-1}h_{1}^{-1}yh_{4}^{-1}y^{-1}h_{1}yh
_{4}$ and $w_{23}=y^{-1}h_{2}^{-1}yh_{3}^{-1}y^{-1}h_{2}yh_{3}$. It is easy to see that $w$ cannot be
reduced to a constant. \parskip0pt

 We claim that
\begin{quote} for any any $g\in F$ satisfying $g(q_{1})<p_{4}$ and $g(p_{4})>q_{1}$ the word $w
_{14}(g)$ is equal to the identity.
\end{quote}
To show this we consider the action of $w_{14}(g)$ on each point from $[0,1]$. Assume, that $t\in
[0,g^{-1}(q_{1}))$. Since $t\notin supp(h_{4})$ we have:
$$w_{14}(g)(t)=g^{-1}h_{1}^{-1}gh_{4}^{-1}g^{-1}h_{1}g(h_{4}(t))=g^{-1}h_{1}^{-1}gh_{4}^{-1}g^{-1}(h_{1}
	(g(t))).$$
By $g^{-1}(h_{1}(g(t)))<g^{-1}(h_{1}(q_{1}))=g^{-1}(q_{1})<p_{4}$ we see $h_{4}^{-1}(g^{-1}(h_{1}(g(t))))
	=g^{-1}(h_{1}(g(t)))$. Thus:
$$w_{14}(g)(t)=g^{-1}h_{1}^{-1}gg^{-1}(h_{1}(g(t)))=t.$$
 If $t\in [g^{-1}(q_{1}),1]$ then since $h_{4}^{\pm 1}(t)\geq min(t,p_{4})$, we have $g(h_{4}^{\pm 1}(t))
	\geq q_{1}$
and hence:
$$w_{14}(g)(t)=g^{-1}h_{1}^{-1}gh_{4}^{-1}g^{-1}h_{1}g(h_{4}(t))=g^{-1}h_{1}^{-1}gh_{4}^{-1}g^{-1}g(h_{4}
	(t))=$$
$$=g^{-1}h_{1}^{-1}g(t)=t.$$
 It follows that for $g$ such that $g(q_{1})<p_{4}$ and $g(p_{4})>q_{1}$, $w_{14}(g)=1_{F}$ and
hence $w(g)=[w_{14}(g),w_{23}(g)]=[1_{F},w_{23}(g)]=1_{F}$. Thus we are left with the case when $g
(q_{1})\geq p_{4}$ (the proof of the case $g(p_{4})\leq q_{1}$ uses the same argument). Now we will
prove that
\begin{quote} for any $g\in F$ satisfying $g(q_{1})\geq p_{4}$ the word $w_{23}(g)$ is equal to the
identity. \end{quote}
 Assume that $t\in [0,g^{-1}(p_{2})]$. Since $g(q_{1})\geq p_{4}$, we have $q_{1}\geq g^{-1}(p_{4})>
g^{-1}(p_{2})\geq t$. Thus:
$$w_{23}(g)(t)=g^{-1}h_{2}^{-1}gh_{3}^{-1}g^{-1}h_{2}gh_{3}(t)=g^{-1}h_{2}^{-1}gh_{3}^{-1}g^{-1}h_{2}g
	(t)=$$
$$=g^{-1}h_{2}^{-1}gh_{3}^{-1}g^{-1}g(t)=t.$$
 Now assume that $t\in (g^{-1}(p_{2}),g^{-1}(q_{2}))$. Then since again $q_{1}\geq g^{-1}(p_{4})> g
^{-1}(q_{2})\geq t$ we obtain:
$$w_{23}(g)(t)=g^{-1}h_{2}^{-1}gh_{3}^{-1}g^{-1}h_{2}gh_{3}(t)=g^{-1}h_{2}^{-1}gh_{3}^{-1}g^{-1}h_{2}
	g(t).$$
Since $h_{2}(g(t))\in (p_{2},q_{2})$ we have $g^{-1}(h_{2}(g(t)))<g^{-1}(q_{2})<q_{1}$ and:
$$w_{23}(g)(t)=g^{-1}h_{2}^{-1}gh_{3}^{-1}(g^{-1}h_{2}g(t))=g^{-1}h_{2}^{-1}g(g^{-1}h_{2}g(t))=t.$$
 Assume that $t\in [g^{-1}(q_{2}),g^{-1}(p_{3})]$. Since we still have $g^{-1}(p_{3})<q_{1}$, we see
that:
$$w_{23}(g)(t)=g^{-1}h_{2}^{-1}gh_{3}^{-1}g^{-1}h_{2}gh_{3}(t)=g^{-1}h_{2}^{-1}gh_{3}^{-1}g^{-1}h_{2}g
	(t)=$$
$$=g^{-1}h_{2}^{-1}gh_{3}^{-1}g^{-1}g(t)=t.$$
 Let $t\in (g^{-1}(p_{3}),g^{-1}(q_{3}))$. Then since $g(p_{3})>q_{2}$ and $h_{3}(t)\neq t\ \Rightarrow\
h_{3}(t)>p_{3}$:
$$w_{23}(g)(t)=g^{-1}h_{2}^{-1}gh_{3}^{-1}g^{-1}h_{2}gh_{3}(t)=g^{-1}h_{2}^{-1}gh_{3}^{-1}g^{-1}gh_{3}
	(t)=$$
$$=g^{-1}h_{2}^{-1}g(t)=t.$$
Finally assume $t\in [g^{-1}(q_{3}),1]$ (and then $g(t)>q_{2}$). Similarly as above:
$$w_{23}(g)(t)=g^{-1}h_{2}^{-1}gh_{3}^{-1}g^{-1}h_{2}gh_{3}(t)=g^{-1}h_{2}^{-1}gh_{3}^{-1}g^{-1}g(h_{3}
	(t))=$$
$$=g^{-1}h_{2}^{-1}g(t)=t.$$
 Now we see that for $g$ such that $g(q_{1})\geq p_{4}$, we have $w_{23}(g)=1_{F}$ and hence
$w(g)=[w_{14}(g),w_{23}(g)]=[w_{14}(g),1_{F}]=1_{F}$. The proof is finished.\\

\ \ \ \ \ \ \ \ \ \ \ \ \ \ \ \ \ \ \ \ \ \ \ \ \ \ \ \ \ \ \ \ \ \ \ \ \ \ \ \ \ \ \ \ \ \ \ \ \ \ \ \ \ \ \ \ \ \ \ \ \ \ \ \ \ \ \ \ \ \ \ \ \ \ \ \ \ $\square$\\

 We now apply the construction from Proposition \ref{lwc} to limits of Thompson's group $F$.

\begin{theo}\label{fre}
	 Suppose we are given a convergent sequence of marked groups $((G_{n}, (x_{0}, x_{1},
	g_{n,1},\ldots , g_{n,s})))_{n<\omega}$, where $G_{n}=F$, $(g_{n,1},\ldots , g_{n,s})\in F$,
	$n\in\mathbb{N}$, and denote by $\mathbb{G}$ its limit. Then $\mathbb{G}\neq F\ast G$
	for any non-trivial $G$.
\end{theo}

 Before the proof we formulate a general statement, which exposes the main point of our argument.

\begin{prop}\label{pro}
	 Let $H=\langle h_{1},\ldots , h_{m}\rangle$ be a finitely generated torsion-free group,
	which satisfies a one variable law with constants and does not satisfy any law without
	constants. Let $\mathbb{G}$ be the limit of a convergent sequence of marked
	groups $((G_{n}, (h_{1},\ldots , h_{m}, g_{n,1},\ldots , g_{n,t})))_{n<\omega}$,
	where $(g_{n,1},\ldots , g_{n,t})\in H$, $G_{n}=H$, $n\in\mathbb{N}$. Then $\mathbb{G}
	\neq H\ast K$ for any non-trivial $K<\mathbb{G}$.
\end{prop}

\emph{Proof.} \ It is clear that $\mathbb{G}$ is torsion-free. To obtain a contradiction suppose that
$\mathbb{G}=H\ast K$, $K\neq\{ 1\}$, and $\mathbb{G}$ is marked by a tuple $(h_{1},\ldots , h_{m},
f_{1},\ldots , f_{t})$. Let $f=u(\bar{h},\bar{f})$ be an element of $K\setminus\{ 1\}$ and let $w(y)$ be
a law with constants in $H$. Obviously $w(u(\bar{h}, g_{n,1},\ldots , g_{n,t}))=1_{H}$ for all $n<
\omega$. It follows from the definition of an $H$-limit group that $w(u(\bar{h},\bar{f}))=1
_{\mathbb{G}}$. Since $w$ was chosen to be non-trivial, with constants from $H$ and $|f|=\infty$,
we obtain a contradiction with the fact that $\mathbb{G}$ is the free product of $H$ and $K$. \\

\ \ \ \ \ \ \ \ \ \ \ \ \ \ \ \ \ \ \ \ \ \ \ \ \ \ \ \ \ \ \ \ \ \ \ \ \ \ \ \ \ \ \ \ \ \ \ \ \ \ \ \ \ \ \ \ \ \ \ \ \ \ \ \ \ \ \ \ \ \ \ \ \ \ \ \ \ $\square$\\

\emph{Proof of Theorem \ref{fre}.} \ It follows directly from Proposition \ref{lwc}, that there is some
word $w(y)$, which is a law with constants in $F$, and hence we just apply Proposition \ref{pro} for $H
=F$, $h_{1}=x_{0}$ and $h_{2}=x_{1}$. \\

\ \ \ \ \ \ \ \ \ \ \ \ \ \ \ \ \ \ \ \ \ \ \ \ \ \ \ \ \ \ \ \ \ \ \ \ \ \ \ \ \ \ \ \ \ \ \ \ \ \ \ \ \ \ \ \ \ \ \ \ \ \ \ \ \ \ \ \ \ \ \ \ \ \ \ \ \ $\square$\\

\section{HNN-extensions}

 Now we proceed to discuss the case of HNN-extensions. For this purpose we consider a sequence
of groups marked by three elements, $(G_{n})_{n<\omega}$, and the corresponding limit group $G =
(\langle x_{0}, x_{1}, g|R_{F} \cup R_{G} \rangle, (x_{0}, x_{1}, g))$. The following theorem is the main
result of the section.

\begin{theo}\label{mt}
	Let $(G_{n})_{n<\omega}$ be a convergent sequence of groups, where $G_{n} =
	(F,(x_{0},x_{1},g_{n}))$, and let $G = (\langle x_{0}, x_{1}, g|R_{F}
	\cup R_{G}\rangle, (x_{0}, x_{1}, g))$ be its limit. Then $G$ is not a
	HNN-extension of Thompson's group $F$ of the following form $\langle x_{0},
	x_{1}, g | R_{F}, ghg^{-1} = h' \rangle$ for some $h,h' \in \emph{F}$.
\end{theo}

In what follows we will need two easy technical lemmas:

\begin{lemm}\label{l1}
	Suppose $g\in F$ and let $x_{0}^{a_{0}}x_{1}^{a_{1}}\ldots x_{n}^{a_{n}}
	x_{n}^{-b_{n}}\ldots x_{1}^{-b_{1}}x_{0}^{-b_{0}}$ be its normal form. There
	is $M\in\mathbb{N}$ such that for all $m>M$:
	$$g^{-1}x_{m}g = x_{m+t} \ \ \emph{or} \ \ gx_{m}g^{-1}=x_{m+t},$$
	where $t = \mid\sum _{i=0} ^{n} (a_{i}-b_{i})\mid$.
\end{lemm}

 \emph{Proof.} Consider the case when $\sum _{i=0} ^{n} (a_{i}-b_{i})\geq 0$. Then for
	sufficently large $m$:

	$$ x_{0}^{b_{0}}x_{1}^{b_{1}}\ldots x_{n}^{b_{n}}x_{n}^{-a_{n}}\ldots x_{1}^{-a_{1}}
	x_{0}^{-a_{0}} x_{m} x_{0}^{a_{0}}x_{1}^{a_{1}}\ldots x_{n}^{a_{n}}x_{n}^{-b_{n}}
	\ldots x_{1}^{-b_{1}}x_{0}^{-b_{0}} =$$

	$$= x_{0}^{b_{0}}x_{1}^{b_{1}}\ldots x_{n}^{b_{n}} x_{m+\sum _{i=0} ^{n} a_{i}}
	x_{n}^{-b_{n}}\ldots x_{1}^{-b_{1}}x_{0}^{-b_{0}} =$$

	$$ = x_{m+\sum _{i=0} ^{n} (a_{i}-b_{i})} $$

	 In the case when $\sum _{i=0} ^{n} (a_{i}-b_{i})<0$ we consider the symmetric
	conjugation and apply the same argument. \\

\ \ \ \ \ \ \ \ \ \ \ \ \ \ \ \ \ \ \ \ \ \ \ \ \ \ \ \ \ \ \ \ \ \ \ \ \ \ \ \ \ \ \ \ \ \ \ \ \ \ \ \ \ \ \ \ \ \ \ \ \ \ \ \ \ \ \ \ \ \ \ \ \ \ \ \ \ $\square$\\

\begin{lemm}\label{l2}
	Under the assumptions of Lemma \ref{l1} the numbers $M$ and $t$ defined in
	that lemma, additionally satisfy the property that for all $m>M$ and $k>0$:
	$$g^{-k}x_{m}g^{k} = x_{m+kt}\ \ \emph{or} \ \ g^{k}x_m g^{-k} = x_{m+kt}.$$
\end{lemm}

\emph{Proof.} If $g^{-1}x_m g=x_{m+t}$ holds (one of possible conclusions of Lemma
	\ref{l1}) then $M\leq m+t$ and applying Lemma \ref{l1} $k$ times we obtain the
	result. The case $gx_m g^{-1} = g_{m+t}$ is similar. \\

\ \ \ \ \ \ \ \ \ \ \ \ \ \ \ \ \ \ \ \ \ \ \ \ \ \ \ \ \ \ \ \ \ \ \ \ \ \ \ \ \ \ \ \ \ \ \ \ \ \ \ \ \ \ \ \ \ \ \ \ \ \ \ \ \ \ \ \ \ \ \ \ \ \ \ \ \ $\square$\\

\emph{Proof of Theorem \ref{mt}.}\ First we prove the theorem in the case of centralized
HNN-extensions. \parskip0pt

 Suppose that $h=h'\neq 1$ in the formulation, i.e. the limit group has a relation of the form $ghg^{-1}
=h$ and denote by $H$ the corresponding HNN-extension of Thompson's group $\langle x_{0},x_{1}, g |
R_{F}, ghg^{-1} = h \rangle$. Assume that $ghg^{-1} = h$ is satisfied in $G$. From the definition of a
limit group it follows, that $g_{n}hg_{n}^{-1}=_{F} h$ for almost all $n$. Denote by $C(h)$ the
centralizer of $h$ and by $C_{1}\oplus\ldots\oplus C_{m}$ its decomposition taken from Fact \ref{GS}.
As almost all $g_{n}$ commute with $h$, almost all $g_{n}$ have a decomposition of the form $g_{n} =
g_{n,1}\ldots g_{n,m}$, where $g_{n,i} \in C_{i}$.\parskip0pt

 As $h\neq 1$, at least one of the factors $C_{1}\oplus\ldots\oplus C_{m}$ is isomorphic to $\mathbb{Z}$,
say $C_{i_{0}}$. Denote by $[a,b]$ the support of elements taken from the subgroup $C_{i_{0}}$. It
follows from the construction of this decomposition, that $h$ can only fix finitely many points in
$[a,b]$.\parskip0pt

 Let us consider the sequence $(g_{n,i_{0}})_{n<\omega}$. Wlog we may suppose that $(g_{n,i_{0}})
_{\omega < \infty}$ consists of powers of some element of $F$ (which is a generator of $C_{i_{0}}$).
Consider the case when it has infinitely many occurences of the same element. If $g'$ occurs infinitely
many times in this sequence, then infinitely many $g_{n}(g')^{-1}$ commute with $x_{[a,b],0}$, $x_{[a,b],
1}$. That gives us a subsequence $(g_{k_{n},i_{0}})_{n<\omega}$ for which the relation $[g_{k_{n}}(g')
^{-1},f]$ holds for all $n$ and for all $f\in \langle x_{[a,b],0},x_{[a,b],1}\rangle$. As $\langle x_{[a,b],0}, x
_{[a,b],1}\rangle$ is isomorphic to Thompson's group $F$, we will find a word of the form $g'y^{-1}f^{-1}
y(g')^{-1}f$ with $f\in \langle x_{[a,b],0}, x_{[a,b],1}\rangle$, which is trivial for $y=\lim _{n\to\infty} g_{k_{n}}$
in the limit group corresponding to the sequence $(g_{k_{n},i_{0}})_{n<\omega}$ and non-trivial for $y=g$
in the group $H$. Indeed, it follows from Britton's Lemma on irreducible words in an HNN-extension
(\cite{LS}, page 181), that the considered word can be reduced in $H$ only if $f^{-1}$ lies in the cyclic
subgroup generated by $h$. But $f\in \langle x_{[a,b],0}, x_{[a,b],1}\rangle$ can be easily choosen outside
$\langle h\rangle$.\parskip0pt

 Let us now assume that the sequence $(g_{n,i_{0}})_{n<\omega}$ is not stabilizing. By the discussion
from the end of Section 1 we see that
\begin{quote} $(\dagger )$\ \ \ \ \ \ \ \ \ \ \ \ \ \ \ for all $m>1$, $[x_{[a,b],0}x_{[a,b],1}^{-1}, x_{[a,b],m}] = 1$.\end{quote}
 On the other hand any $g_{n,i_{0}}$ is a power of some fixed element from $C_{i_{0}}$. Thus we see by
Lemma \ref{l2} that for $M$ found for the generator of $C_{i_{0}}$ as in Lemma \ref{l1}:
	$$(\forall n)\ g_{n,i_{0}}^{-1}x_{[a,b],M}g_{n,i_{0}}=x_{[a,b],M+t_{n}}\ \ \emph{or}\ \ (\forall n)\
		g_{n,i_{0}}x_{[a,b],M}g_{n,i_{0}}^{-1}=x_{[a,b],M+t_{n}}, \ t_{n}\geq 0.$$
 Substituting $x_{[a,b],M+t_{n}}$ into $(\dagger )$ instead of $x_{[a,b],m}$ we have:
	$$[x_{[a,b],0}x_{[a,b],1}^{-1}, g_{n,i_{0}}^{-1}x_{[a,b],M}g_{n,i_{0}}] = 1\ \ \emph{or}\ \ [x_{[a,b],0}x
		_{[a,b],1}^{-1}, g_{n,i_{0}}x_{[a,b],M}g_{n,i_{0}}^{-1}] = 1.$$
 Thus we see that one of the following relations holds for all $n$:
	$$ [x_{[a,b],0}x_{[a,b],1}^{-1}, g_{n}^{-1}x_{[a,b],M}g_{n}] = 1 \ \ \emph{or} \ \ [x_{[a,b],0}
		x_{[a,b],1}^{-1}, g_{n}x_{[a,b],M}g_{n}^{-1}] = 1.$$
 Suppose that the first relation holds for all $n$'s.
Consider the corresponding word in the group $H$:
	$$ w = x_{[a,b],1}x_{[a,b],0}^{-1}g^{-1}x_{[a,b],M}^{-1}gx_{[a,b],0}x_{[a,b],1}^{-1}g^{-1}x_{[a,b],
		M}g.$$
 We claim that $w \neq 1$ in this HNN-extension. Once again, it follows from Britton's Lemma on irreducible
words in an HNN-extension, that we can reduce $w$ if $x_{[a,b],M}$ is a power of $h$ or $x_{[a,b],0}
x_{[a,b],1}^{-1}$ is a power of $h$. We know that $x_{[a,b],0},	x_{[a,b],1}, \ldots x_{[a,b],m}, \ldots$
generate the group $F _{[a,b]}$, which is isomorphic to $F$. From the properties of $F$ we know that for
different $m,m'>M$, $x_{[a,b],m}$ and $x_{[a,b],m'}$ do not have a common root. Thus, possibly increasing
the number $M$, we can assume that $x_{[a,b],M}$ is not a power of $h$. If $x_{[a,b],0}x_{[a,b],1}^{-1} =
h^{d}$ for some integer $d$, then $h^{d}$ fixes pointwise the segment $[\frac{1}{4}a + \frac{3}{4}b,b]
\subset [a,b]$. Hence $h$ also fixes some final subinterval of [a,b]. This gives a contradiction as $h$ was
chosen to fix only finitely many points in $[a,b]$. This finishes the case of centralized
HNN-extensions.\parskip0pt

 Generally, let us consider the situation, where in the limit group we have one relation of the form
$ghg^{-1}=h'$ for some $h, h'\in F$. By the construction of limit groups $h'=h^{f}$ for some element $f
\in F$. Indeed, if $h$ and $h'$ are not conjugated in $F$, then there is no sequence $(g_{n})_{n<\omega}$
in $F$ with $g_{n}hg_{n}^{-1}=h'$ for almost all $n$. We now apply the argument above: let $(fg_{n})
_{n<\omega}$ be a sequence convergent to the element $fg$. It obviously commutes with $h$, so we can
repeat step by step the proof above. That completes the proof. \\

\ \ \ \ \ \ \ \ \ \ \ \ \ \ \ \ \ \ \ \ \ \ \ \ \ \ \ \ \ \ \ \ \ \ \ \ \ \ \ \ \ \ \ \ \ \ \ \ \ \ \ \ \ \ \ \ \ \ \ \ \ \ \ \ \ \ \ \ \ \ \ \ \ \ \ \ \ $\square$\\

\end{document}